# Independent graph of the finite group

T. Chalapathi [1]     R. V M S S Kiran Kumar [2,*]

[1] Assistant Professor, Department of Mathematics, Sree Vidyanikethan Eng.College Tirupati,-517502, Andhra Pradesh, India.

[2] Research Scholar, Department of Mathematics, S.V.University, Tirupati,-517502, Andhra Pradesh, India

*Corresponding author* Email: [b]kksaisiva@gmail.com.

**Abstract**

Let $a$ and $b$ be any two elements in the group $Z_n$ of integers modulo $n$. Then $a$ and $b$ are called independent if $o(a) \neq o(b)$. In this paper, we introduce and study independent graph of the group $Z_n$, denoted by $I_G(Z_n)$, is undirected simple graph whose vertex set is $Z_n$ and two distinct vertices $a$ and $b$ are adjacent in $I_G(Z_n)$ if and only if $a$ and $b$ are independent in $Z_n$.

1. Introduction

The multidisciplinary research between finite algebraic structures and graphs has been the most productive area of algebraic graph theory. Algebraic graph theory is a branch of modern mathematics in which algebraic methods are applied to real world problems about graphs. These graphs are a nice composition of three main branches of mathematics, viz., number theory, abstract algebra and graph theory. The author (Frucht, 1994) shown that all finite groups can be represented as the automorphism group of a connected algebraic graph. In 1878, the author

Cayley was originated an idea between finite groups and graphs. According to the author Cayley, a graphical representation of a finite group is given by a set of generators and their relations. It provided an algebraic method of visualizing a finite group and connects many important branches of mathematics. For algebraic graphs reader may refer, (Godsil & Royle, 2001).

The theory of finite groups occupies a central position in mathematics for studying symmetries of the objects in the real world. Now a day, finite groups play a specific role in different areas such as algebraic coding theory, algebraic cryptography, design theory, and engineering science. In fact, author Gauss introduced the finite cyclic group $Z_n$ of integers 0, 1, 2 ,..., $n-1$ with respect to addition modulo $n$. For further details of the group $Z_n$ reader refer (Walter, 1996a; Jungnickel, 1996b; Alkam & Abu Osba, 2008a; & Chalapathi & Kiran Kumar 2015). For given a positive integer $n$ and an element $a$ of $Z_n$ such that $0 \leq a < n$. The order of $a$ denoted by $o(a)$ and defined as a least positive integer $k$ less than or equal to $n$ such that $ka = 0$. Let the group $Z_n$ contains the subsets $U_n$, $S_n$ and $N_n$, where $U_n$ is the set of all group units, $S_n$ is the set of additive involutions of $Z_n$ and $N_n = Z_n - (U_n \cup S_n)$. Note that $N_n = \phi$ if and only if $n$ is prime and $|U_n| = \varphi(n)$, the Euler-totient function of $n$. Here $U_n \cap S_n = \phi$, $S_n \cap N_n = \phi$, $N_n \cap U_n = \phi$ and $-U_n = U_n$, where $-U_n = \{-u : u \in U_n\}$. For instance, $U_6 = \{1, 5\}$, $S_6 = \{0, 3\}$ and $N_6 = \{2, 4\}$.

The main motivation of our paper is (Carrie Finch & Lenny Jones, 2002; Herzog, 1977; Freud & Pál Pálfy, 1996c; Farrokhi & Saeedi, 2014). In (Carrie Finch & Lenny Jones, 2002), the authors Finch and Jones introduced a nice connection between Fermat number and finite groups. Also, they defined order subset (OS) and perfect order subset (POS) of finite groups. The OS of a finite group $G$ determined by an element $x \in G$ is defined to be the non-empty set $OS(x) = \{y \in G : o(y) = o(x)\}$ and similarly $G$ is said to be a POS- group if for each $x \in G$ the cardinality of $OS(x)$ is a devisor of the order of $G$. For further investigations of POS-groups reader refer (Kumar Das, 2009a). In (Freud & Pál Pálfy, 1996c), the authors Róbert Freud and Péter Pál Pálfy introduced the relation $f(k, G) = \varphi(k)s(k, G)$, where $f(k, G)$ denote number of elements of order $k$, $s(k, G)$ denote the number of cyclic subgroups of a finite group $G$ and $\varphi(k)$ is the Euler- totient function of $k$.

In finite group theory, many researchers study and analyze different finite cyclic groups and their structural properties while in algebraic graph theory mainly we focus on the group theoretic graphs. Inspired by (Chalapathi & Kiran, 2015), Dieter Jungnickel studied so many new properties of orders elements of finite cyclic group. Further, the study of algebraic graphs using the properties of finite cyclic groups and their graphs has become an inspiring research in the recent years. In this sequel, we study an undirected simple graph related to finite cyclic group $Z_n$. We call this graph is an independent graph, denoted by $I_G(Z_n)$ with vertex set $Z_n$ and two distinct vertices $a$ and $b$ in $Z_n$ are adjacent in $I_G(Z_n)$ if and only if $a$ and $b$ are independent in $Z_n$.

### 2. Properties of the independent graph

In this section, we begin definition of independent graph and establish the nature of the degree of each vertex, and we obtain a formula for enumerating size of independent graphs. Also, we characterize $n$ for which the graph $I_G(Z_n)$ is complete, star and Hamiltonian. Further, we compute girth and diameter of $I_G(Z_n)$.

**Definition. 2.1** Let $n > 1$ be a positive integer. Then the undirected simple graph $I_G(Z_n)$ is called independent graph with vertex set $Z_n$ and two distinct vertices $a$ and $b$ are adjacent in $I_G(Z_n)$ if and only if $o(a) \neq o(b)$, where $o(a), o(b)$ denotes orders of $a$ and $b$ respectively.

**Example.2.2** The following graph is the independent graph $I_G(Z_6)$ of the group $Z_6$

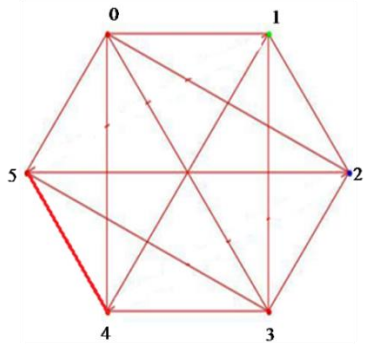

Fig. 2.3 Independent graph $I_G(Z_6)$

We start with some graph theoretical properties of $I_G(Z_n)$ for each value of $n > 1$.

**Theorem.2.4** For any positive integer $n > 1$, the independent graph $I_G(Z_n)$ is connected.

**Proof.** It is obvious since $o(0) = 1$ and $o(a) \neq 1$ for every $a \neq 0$ in $Z_n$, so the vertex 0 is independent with remaining all the vertices of the graph $I_G(Z_n)$, and hence $I_G(Z_n)$ is connected.

Before going to further results of $I_G(Z_n)$ we give the following basic lemmas which will be give an enumeration of total number of elements in $S_n$ and $N_n$ for each $n > 1$.

**Lemma. 2.5** The total number of additive involutions of the group $Z_n$ is given by

$$|S_n| = \begin{cases} 1 \text{ if } n \text{ is odd} \\ 2 \text{ if } n \text{ is even} \end{cases}.$$

**Proof.** Let $n > 1$ be a positive integer. Then we consider the following two cases on $n$.

**Case. (i)** Let $n$ be odd. Suppose $a \neq 0$ be an element in $S_n$. Then clearly, $a \notin U_n$ and $a \notin N_n$. If $a \in U_n$ then there exist $a^{-1}$ in $U_n$ such that $a^{-1}(2a) = 0$. This implies $2 = 0$, that is $n$ must be even, which is a contradiction to our hypothesis that $n$ is odd. Similarly, if $a \in N_n$, then there exist again a contradiction on $n$. So, we have $2a = 0$ implies that $a = n - a$. This shows that $a = \dfrac{n}{2}$ and thus $n$ must be even, which is also not possible because $n$ is odd. So our assumption that $a \neq 0$ is not true. Hence $a = 0$ is the only one element in $S_n$, and $|S_n| = 1$ when $n$ is odd.

**Case. (ii)** Let $n$ be even. Trivially, $\dfrac{n}{2} + \dfrac{n}{2} \equiv 0 \pmod{n}$. Therefore, by the definition of set of additive involutions,

$$S_n = \{s \in Z_n : 2s = 0\}$$

$$= \{s \in Z_n : 2s \equiv n \pmod{n}\}$$

$$= \left\{s \in Z_n : s = 0, \dfrac{n}{2}\right\} = \left\{0, \dfrac{n}{2}\right\}. \text{ Hence, } |S_n| = 2.$$

**Lemma. 2.6** Let $n > 1$ be a positive integer. Then $|N_n| = \begin{cases} n - \varphi(n) - 1 \text{ if } n \text{ is even} \\ n - \varphi(n) - 2 \text{ if } n \text{ is odd} \end{cases}.$

**Proof.** For any positive integer $n > 1$, we have

$$Z_n = \{0, 1, 2,..., n-1\}, \ U_n = \{u \in Z_n : gcd(u, n) = 1\},$$

$$S_n = \{s \in Z_n : s = -s\} \text{ and } N_n = \{a \in Z_n : \gcd(a, n) \neq 1 \text{ and } a \neq -a\}.$$

Clearly, $U_n \cap S_n = \phi$, $S_n \cap N_n = \phi$ and $N_n \cap U_n = \phi$. Therefore, $|Z_n| = |U_n| + |S_n| + |N_n|$. In view of Lemma [2.5] and $|U_n| = \varphi(n)$ follows the result.

**Theorem. 2.7** If $a$ is any vertex in $I_G(Z_n)$, then the degree of $a$ is given by

$$\deg(a) = \begin{cases} n-1 & \text{if } a \in S_n \\ n - \varphi(n) & \text{if } a \in U_n \\ \varphi(n) + 2 \text{ or } \varphi(n) + 1 & \text{if } a \in N_n \end{cases}.$$

**Proof.** Suppose $n > 1$, and let $a$ be a vertex of the graph $I_G(Z_n)$. Since, $Z_n = U_n \cup S_n \cup N_n$. If $a \in S_n$, then $2a = 0$ imply $a = -a$. Therefore, $o(a) = o(-a)$. It is clear that the vertex $a$ in $S_n$ is adjacent to all other vertices $b$ which is not in $S_n$, because $o(a) \neq o(b)$ and $1 \leq o(a) \leq 2$ for every $a \in S_n$ and $b \notin S_n$. This implies that the degree of $a$ is $n-1$. If $a \in U_n$ then there exist $a^{-1}$ in $U_n$ such that $aa^{-1} = 1 = a^{-1}a$, so that $o(a) = o(a^{-1}) = o(n-a) = n$, since $\gcd(a, n) = 1$ and $\gcd(n-a, n) = 1$. This implies that $a$ is not adjacent to $a^{-1}$ in $I_G(Z_n)$, which are $\varphi(n)$ in total. Hence, the degree of $a$ is $n - \varphi(n)$. Finally, if $a \in N_n$ and $N_n \neq \phi$ then $a \notin U_n \cup S_n$, it is clear that $-a \in N_n$. This implies that $1 < o(a) < n$. So, there exist $u \in U_n$ such that $o(a) \neq o(u)$, where $o(u) = n$. Also, if $s \in S_n$, then $2s = 0$ but $2a \neq 0$. It is clear that $o(a) \neq o(s)$. Therefore, each pair vertices in $N_n = Z_n - (U_n \cup S_n)$ are not adjacent because $o(a) = o(b)$ for every $a, b \in N_n$. So, by the Lemma [2.6], degree of $a$ is $n - |N_n| = \varphi(n) + 2$ or $\varphi(n) + 1$. Hence the result follows.

Now we recall two fundamental results from (Jungnickel, 2009b) for simple graphs, and we shall investigate some important concrete properties of $I_G(Z_n)$.

**Theorem. 2.8** Let $\deg(a_i)$ denote degree of vertex $a_i$ in the graph $G$. Then $2|E| = \sum_{i=1}^{n} \deg(a_i)$.

**Theorem. 2.9** Let $G$ be a complete graph of order $n$. Then the size of $G$ is $\binom{n}{2} = \frac{n(n-1)}{2}$.

**Theorem. 2.10** Let $n > 1$ be a positive integer. Then the number of edges in the independent graph

$I_G(Z_n)$ is $|E(I_G(Z_n))| = \begin{cases} \frac{1}{2}\left((n-1)^2 - \varphi(n)(\varphi(n)-2)\right) & \text{if } n \text{ is odd} \\ \frac{1}{2}\left((n-1)^2 + 1 - \varphi(n)(\varphi(n)-2)\right) & \text{if } n \text{ is even} \end{cases}$.

**Proof.** Write $N_n = Z_n - U_n \cup S_n$. For each $n > 1$ we have the group $Z_n$ can be written as disjoint union of $U_n$, $S_n$ and $N_n$. By the Theorem [2.8],

$$2|E(I_G(Z_n))| = \sum_{s \in S_n} \deg(s) + \sum_{u \in U_n} \deg(u) + \sum_{a \in N_n} \deg(a).$$

First we consider $n$ is odd. By the Lemma [2.5] and Lemma [2.6], we have $|S_n| = 1$, $|N_n| = n - \varphi(n) - 1$. But, $|U_n| = \varphi(n)$. Thus, in view of Theorem [2.8],

$$|E(I_G(Z_n))| = \frac{1}{2}\left(1(n-1) + \varphi(n)(n-\varphi(n)) + (n-\varphi(n)-1)(n-2)\right)$$

$$= \frac{1}{2}\left((n-1)^2 - \varphi(n)(\varphi(n)-2)\right).$$

Next we consider $n$ is even. Again by the Lemma [2.5], $|S_n| = 2$. Also, $|U_n| = \varphi(n)$ and, by the Lemma [2.6], $|N_n| = n - \varphi(n) - 2$. Therefore,

$$|E(I_G(Z_n))| = \frac{1}{2}\left(2(n-1) + \varphi(n)(n - \varphi(n)) + (n - \varphi(n) - 2)(n - 2)\right)$$

$$= \frac{1}{2}\left((n-1)^2 + 1 - \varphi(n)(\varphi(n) - 2)\right).$$

**Example. 2.11** The size of $I_G(Z_4)$ and $I_G(Z_5)$ is 5 and 4 respectively.

**Theorem. 2.12** The independent graph $I_G(Z_n)$, $n > 2$ is never complete.

**Proof.** Suppose on contrary that $I_G(Z_n)$, $n > 2$ is a complete graph. Then, by the Theorem [2.9], the total number of edges in a simple graph of order $n$ is $\frac{n(n-1)}{2}$, but in view of Theorem [2.10], we arrived a contradiction to the compactness of $I_G(Z_n)$.

**Corollary. 2.13** The independent graph $I_G(Z_n)$ is complete if and only if $n = 2$.

**Proof.** It is clear from the fact that $o(0) = 1$ and $o(1) = 2$ in the group $Z_2$.

In view of Theorem [2.12], the following result is obvious.

**Theorem. 2.14** Let $p$ be a prime. Then the $I_G(Z_p)$ is a star graph, which is isomorphic to $K_{1, p-1}$.

**Theorem. 2.15** Let $n > 1$ be a positive integer. Then $gir(I_G(Z_n)) \in \{3, \infty\}$.

**Proof.** First suppose that $n$ is prime. Then, by Theorem [2.14], graph $I_G(Z_n)$ is isomorphic to $K_{1,n-1}$, it is acyclic graph, and hence $gir(I_G(Z_n)) = \infty$. Now suppose that $n$ is composite. Then there exist a proper divisor $m$ of $n$ such that $1 < m < n$. It is easy to see that $o(0) \neq o(a)$, $o(a) \neq o(b)$ and $o(b) \neq o(0)$ when $o(0) = 1$, $o(a) = m$ and $o(b) = n$. Hence, $0-a-b-0$ is a three cycle, which is smallest in $I_G(Z_n)$. Further suppose $a$ and $b$ are non- identity elements in $I_G(Z_n)$, $n > 4$, there exist another vertex $a+b$ in $I_G(Z_n)$ such that $o(a) \neq o(a+b)$ and $o(a+b) \neq o(b)$ because $o(a+b) | Lcm[o(a), o(b)]$. Thus, $a-(a+b)-b-a$ is a there cycle in $I_G(Z_n)$. Hence, girth of $I_G(Z_n)$ is 3.

For distinct vertices $a$ and $b$ of a simple graph $G$, the diameter of $G$ defined by

$$\operatorname{diam}(G) = \max \{d(a, b) : a, b \in V(G)\},$$

where $d(a, b)$ is the length of the shortest path between $a$ and $b$.

**Theorem.2.16** The diameter of an independent graph, $I_G(Z_n), n > 1$ is at most 2.

**Proof.** We know that the independent graph $I_G(Z_n)$, $n > 1$ having the vertices of the form $0, 1, 2, ..., n-1$. But the vertex 0 is adjacent to every vertex of $I_G(Z_n)$ so that there exist a path between the vertices 0 and $a \neq 0$ in $I_G(Z_n)$, and thus $d(0, a) \geq 1$. Now, suppose $a \neq 0$ and $b \neq 0$ be any two vertices in $I_G(Z_n)$. If $a$ is adjacent to $b$, then, obviously $d(a, b) = 1$. However, $a$ is not adjacent to $b$ for all $0 < a, b \leq n-1$, so, $d(a, b) > 1$, but in $I_G(Z_n)$, there always exists a shortest path $a-0-b$ of length 2, which gives $d(a, b) = 2$ for every two non-adjacent vertices $a \neq 0$ and $b \neq 0$ in $I_G(Z_n)$. It follows that diameter of $I_G(Z_n)$ is 2. Hence, $\operatorname{diam}(I_G(Z_n)) \leq 2$.

A cycle in a simple graph $G$ is called Hamiltonian cycle if it visits every vertex exactly once, and such a graph $G$ is called Hamiltonian. In the graph theory to find the Hamiltonian cycle is an NP-complete problem. For more details on Hamiltonian graph and its cycles the reader is referred to (Bondy & Murty, 2008b). In fact, one of the interesting properties of independent graph is that they provide a class of Hamiltonian graph.

**Theorem. 2.17** Let $n \geq 4$ be a composite number. Then $I_G(Z_n)$ is Hamiltonian.

**Proof.** For each composite number $n \geq 4$, we have to show that the graph $I_G(Z_n)$ is Hamiltonian. For this we shall show that $I_G(Z_n)$ properly contains a cycle of length $n$. To do this, let $i$ and $i+1$ be any two consecutive vertices in $I_G(Z_n)$, then $o(i) \neq o(i+1)$ for each $i$, $0 \leq i \leq n-1$. Suppose, $o(i) = o(i+1)$. Then there exist at least positive integer $k$ such that $ki = k(i+1)$. This implies that $k = 0$, it is not possible because $k > 0$. Thus there exist an edge between two vertices $i$ and $i+1$ in $I_G(Z_n)$ for each $0 \leq i \leq n-1$. So, we must construct a cycle $0-1-2-\cdots-(n-2)-(n-1)-0$ which covers all the vertices in $I_G(Z_n)$, and thus it is a Hamilton cycle of length $n$ in $I_G(Z_n)$. Hence $I_G(Z_n)$ is Hamiltonian.

The following remark is obvious from the definition of $I_G(Z_n)$.

**Remark. 2.18** $I_G(Z_n)$ is not Hamiltonian if and only if $n \in \{2, 3\}$.

3. **The independent graphs are partite**

In this section, we show that the independent graphs are *bipartite* and *complete d(n)-partite* graphs for various values of $n > 1$.

In the field of graph theory, a bipartite graph $G$ is a simple undirected graph whose vertex set $V$ can be divided into two independent sets $A$ and $B$ such that every edge in $A$ connects to every edge in $B$ only. Hence, $A$ and $B$ together form a bipartition $(A, B)$ of $G$. This bipartition of a graph implies the validity of following results.

**Theorem. 3.1** For each prime $p$, the independent graph $I_G(Z_p)$ is bipartite.

**Proof.** Write $A = \{0\}$ and $B = \{a \in Z_p : o(a) = p, a \neq 0\}$. Then $A$ and $B$ are non empty subsets of the vertex set of $I_G(Z_p)$ such that $A \cup B = Z_p$ and $A \cap B = \phi$. Because of $o(0) \neq o(a)$, for every vertex $a \neq 0$ in $I_G(Z_p)$, the pair $(A, B)$ is bipartition of the graph $I_G(Z_p)$, and hence $I_G(Z_p)$ is a bipartite graph.

**Theorem. 3.2** If $n$ is a composite number, then $I_G(Z_n)$ is not a bipartite graph.

**Proof.** Assume that $n$ is composite. Write $U_n' = Z_n - U_n$, where $U_n$ be the set of group units of the additive group $Z_n$. We now show that $I_G(Z_n)$ is not a bipartite graph. Suppose, $I_G(Z_n)$ is bipartite. Then there exist a bipartition $(U_n, U_n')$ in $I_G(Z_n)$. Obviously, $1 \in U_n$ and $0 \in U_n'$. Since $n$ is a composite number and $Z_n$ is a cyclic group of order $n$, so there exist at least one proper divisor $d$ of $n$ such that, by the Lagrange's theorem of groups (Lanski, 2010), $o(0) \mid d$ and $d \mid o(1)$. Therefore, $d$ can be placed neither in $U_n$ nor in $U_n'$. This violates the condition of the bipartition of the bipartite graph. Hence, $I_G(Z_n)$ is not a bipartite graph for any composite number $n$.

**Theorem.3.3** Let $d(n)$ be the total number of divisors of $n$. Then the independent graph $I_G(Z_n)$ is a complete $d(n)$ – partite graph.

**Proof.** Set $D(n) = \{d : d \mid n\}$. Then $|D(n)| = d(n)$, total number of divisors of $n$. For each divisor $d_i$, we considered the set $D_i = \{a \in Z_n : o(a) = d_i\}$ where $1 \leq i \leq d(n)$. So, the group $Z_n$ is the disjoint union of $D_1, D_2, \ldots, D_{d(n)}$, that is, $Z_n = D_1 \cup D_2 \cup \ldots \cup D_{d(n)}$. It is clear that each $D_i$ is independent vertex subset of $I_G(Z_n)$. But the vertex $a$ is not adjacent to the vertex $b$ in $I_G(Z_n)$ if and only if $a$ and $b$ are both lies in $D_i$ because $a, b \in D_i$ if and only if $o(a) = o(b)$. Now suppose $a \in D_i$ and $b \in D_j$ for $i \neq j$. Then clearly, $o(a) = d_i$ and $o(b) = d_j$. If either $i < j$ or $i > j$, then $d_i \neq d_j$, that is, $o(a) \neq o(b)$. Therefore, $a$ is adjacent to $b$ in $I_G(Z_n)$ and thus the edge joining two independent sets of vertices means that each vertex is one independent set $D_i$ is adjacent to each vertex of the other independent set $D_j$ in $I_G(Z_n)$. Note that the number of independent sets in $I_G(Z_n)$ is $d(n)$. Hence $I_G(Z_n)$ is a complete $d(n)$ – partite graph.

**Corollary.3.4** Let $p$ and $q$ be two distinct primes. Then the independent graph $I_G(Z_{pq})$ is a complete $4$ – partite graph.

**Proof.** Write $D_1 = \{a \in Z_{pq} : o(a) = 1\}$, $D_p = \{a \in Z_{pq} : o(a) = p\}$, $D_q = \{a \in Z_{pq} : o(a) = q\}$ and $D_{pq} = \{a \in Z_{pq} : o(a) = pq\}$. Because $1 \neq p$, $p \neq q$ $q \neq pq$ and $pq \neq 1$, the sets $D_1, D_p, D_q$ and $D_{pq}$ are independent vertex subsets of $I_G(Z_{pq})$. Therefore, the pairs $(D_1, D_p)$, $(D_p, D_q)$ $(D_q, D_{pq})$ and $(D_{pq}, D_1)$ are bipartitions of the graph $I_G(Z_{pq})$, and hence $I_G(Z_{pq})$ is a complete $4$ – partite graph.

**Example.3.5** If an edge between two independent sets $D_i$ and $D_j$ means that each vertex in $D_i$ is adjacent to each vertex in $D_j$ for $i \neq j$. Then independent graph $I_G(Z_6)$ is a complete $4$-partite graph which is shown below.

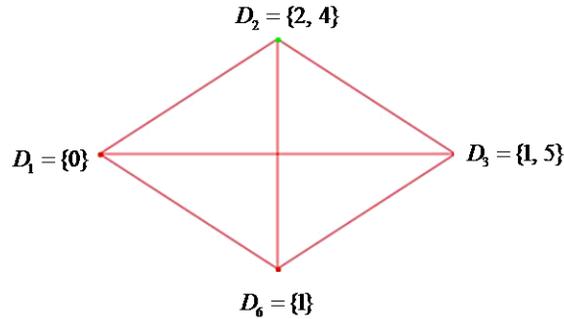

4. **Clique and chromatic numbers of $I_G(Z_n)$**

In this section, we determine the clique and chromatic numbers of $I_G(Z_n)$, and hence to prove that $I_G(Z_n)$ is strongly perfect graph.

First we note that the following. The clique of a simple graph $G$ is a maximal complete subgraph, and clique number $\omega(G)$ of $G$ is the total number of mutually adjacent vertices in $G$. An independent set of vertices in $G$ is the set of pair wise non-adjacent vertices, and the independent set of $G$ is also called co-clique. The chromatic number of colors required to color the vertices of $G$, that is the chromatic number of $G$ is $\chi(G) = \min\{l : G \text{ is } l-\text{colorable}\}$. If the chromatic number is equal to clique number of $G$, then $G$ is called weakly perfect. Otherwise $G$ is called strongly perfect.

**Theorem.4.1** For any positive integer $n > 2$, the clique number $\omega(I_G(Z_n))$ of $I_G(Z_n)$ is

$$\frac{1}{2}\left(n - \varphi(n)(\varphi(n) - 2) + |S_n|\right)$$

**Proof.** Let $n > 2$ be a positive integer. Then each vertex in $S_n$ is adjacent with all other vertices of the graph $I_G(Z_n)$, because $1 \leq o(s) \leq 2$ for every $s \in S_n$. Now for any two vertices $i$ and $j$ in $I_G(Z_n)$, we consider the following two cases on $i$ and $j$.

**Case (i).** Suppose each vertex $i$ is non-adjacent with exactly one vertex $j$. Then $i$ and $j$ must be in $N_n = Z_n - U_n \cup S_n$, by the Lemma [2.6], such vertices are $|N_n| = n - \varphi(n) - |S_n|$. Clearly, in this case, the pair of non-adjacent vertices in $I_G(Z_n)$ are $\frac{1}{2}(n - \varphi(n) - |S_n|)$.

**Case (ii).** Suppose each vertex $i$ is non-adjacent with more than one vertex $j$. Then $i$ and $j$ must be in $U_n$, such vertices are $|U_n|(|U_n| - 1)$, that is, $\varphi(n)(\varphi(n) - 1)$ because $|i| = |j| = n$ if and only if $i, j \in U_n$. So, in this case the pair of non-adjacent vertices are $\frac{1}{2}\varphi(n)(\varphi(n) - 1)$.

Clearly, from the above two cases, the pair of non-adjacent vertices are $\frac{1}{2}(n - \varphi(n) - |S_n|) + \frac{1}{2}\varphi(n)(\varphi(n) - 1)$, and hence total number of mutually adjacent vertices in $I_G(Z_n)$ are $n - \frac{1}{2}(n - \varphi(n) - |S_n|) - \frac{1}{2}\varphi(n)(\varphi(n) - 1) = \frac{1}{2}(n - \varphi(n)(\varphi(n) - 2) + |S_n|)$, which is a size of maximum clique.

**Remark. 4.2** From (Bondy & Murty, 2008b), we have $\chi(G) \geq \omega(G)$. Using together with the Theorem [4.1], we obtain the inequality $\chi(I_G(Z_n)) \geq \frac{1}{2}(n - \varphi(n)(\varphi(n) - 2) + |S_n|)$.

**Theorem.4.3** Let $n > 2$ be a positive integer then the chromatic number $\chi(I_G(Z_n))$ of the graph $I_G(Z_n)$ is $\frac{1}{2}(3n - 3\varphi(n) + |S_n|)$.

**Proof.** Except the vertices in $U_n \cup S_n$, every vertex in the graph $I_G(Z_n)$ is non-adjacent with precisely one vertex, that is its inverse, because $o(a) = o(-a)$ for every $a \in N_n$. Therefore, maximum independent set in $I_G(Z_n)$ is of size 2, moreover such independent sets are $\frac{1}{2}(n - \varphi(n) - |S_n|)$ in number. However, each independent set is uniquely colorable, which means that for all these vertices we need $\frac{1}{2}(n - \varphi(n) - |S_n|)$ colors. Further, vertices in $S_n$ are adjacent with all other vertices, thus we require $|S_n|$ more colors distinct from $\frac{1}{2}(n - \varphi(n) - |S_n|)$ colors. Similarly, vertices in $U_n$ are also adjacent with all the vertices in $Z_n - U_n$, because $o(u) \neq o(a)$ where $o(u) = n$ for every $u \in U_n$ and $1 \leq o(a) < n$ for every $a \in Z_n - U_n$, and thus we further require $|Z_n| - |U_n|$, that is, $n - \varphi(n)$ more colors. Hence, minimum number of colors required to color the independent graph $I_G(Z_n)$ are $\frac{1}{2}(n - \varphi(n) - |S_n|) + |S_n| + n - \varphi(n)$, which is equal to $\frac{1}{2}(3n - 3\varphi(n) + |S_n|)$.

**Theorem. 4.4** For each $n > 2$, the graph $I_G(Z_n)$ is strongly perfect.

**Proof.** In view of Theorems [4.2] and [4.3], result follows.

**Acknowledgments**

The authors express their sincere thanks to Prof.L.Nagamuni Reddy and Prof.S.Vijaya Kumar Varma for his suggestions during the preparation of this paper and the referee for his suggestions.